\documentclass[11pt]{amsart}

\usepackage[margin=1in]{geometry}
\usepackage{amsmath,amssymb}
\usepackage{hyperref}
\usepackage{setspace}

\title{\textbf{Diamond Open Access: The AMR Experiment}}
\author{Alex Kontorovich}
\address{Department of Mathematics, Rutgers University}
\email{alex.kontorovich@rutgers.edu}
\thanks{Partially supported by NSF grant DMS-2302641.}

\begin{document}

\maketitle

\begin{abstract}
    Diamond open access journals charge neither readers nor authors. Despite long-standing support for this ideal within mathematics, relatively few such journals exist. This article documents the Association for Mathematical Research’s experience building and operating diamond open access journals, focusing on the infrastructure, cost, and editorial practices that make the model viable. It aims to clarify why earlier reform efforts have been difficult to replicate and how a lightweight institutional framework can lower the barrier to adoption.
\end{abstract}

\section{Introduction}

I had been thinking about journal publishing for some time. Like many mathematicians, I signed
Tim Gowers’s 2012 call \cite{CostOfKnow2012} to boycott Elsevier.
Years later, however, little had changed. In July 2017, Gowers wrote a blog post \cite{GowersBlog2017}
celebrating that the editorial board of the \emph{Journal of Algebraic Combinatorics} was
``flipping''—resigning en masse from Springer to start a new diamond open access journal,
\emph{Algebraic Combinatorics}, under the same editorial leadership. Yet even as he celebrated,
Gowers admitted: ``Although it is very exciting that one journal is flipping, I must also admit to
disappointment at how low our strike rate has been so far.''

One of the most serious attempts to reform mathematical publishing over the past several decades
was Rob Kirby’s founding of Mathematical Sciences Publishers (MSP) in spring 2004 (see \cite{Kirby2004}). MSP became
the home of \emph{Geometry \& Topology}, a high-quality journal created in 1996 explicitly in response to rising
subscription prices. MSP’s journals are nonprofit and fairly priced, dramatically less expensive than those of commercial publishers; that said, they remain
subscription-based and therefore not fully open.

My own experience brought the issue into sharper focus. When I was editor-in-chief of
\emph{Experimental Mathematics}, I once tried to access a paper
from my own journal and was asked to pay \$40. Around the same time, I had published a paper in
Gowers’s diamond open access journal \emph{Discrete Analysis}, and saw how that model worked
beautifully: no fees to authors, no fees to readers, and a fully functional journal. This raised a
simple question: if this works so well, why hasn’t it spread further?

The answer, I eventually came to realize, is not a lack of goodwill or conviction. It is a coordination problem.
Creating a ``real'' journal from scratch—even one that exists only online—requires solving many
small but nontrivial problems: DOIs, ISSNs, editorial management systems, hosting, indexing,
archiving, and long-term continuity. Each component is manageable in isolation, but assembling
them all without an existing template is daunting enough to deter most editorial boards from
trying.

In the fall of 2021, the fledgling Association for Mathematical Research\footnote{\url{https://amathr.org}} (AMR) contacted me. Its
founders were interested in whether the technologies mathematicians had grown accustomed to
during the pandemic might change how mathematics is practiced and communicated. This
provided a golden opportunity to attack the problem I had been thinking about for so long.

Over the next year, I interviewed about a hundred mathematicians around the world to understand whether this was a problem people genuinely wanted solved, and what it would take to solve it. The result was the AMR’s first journal, \emph{JAMR} (\url{https://jamathr.org}), and more importantly, the development of a \textbf{replicable process}. After three years of operation, \emph{JAMR} is now indexed by Scopus and reviewed by MathSciNet, reflecting its status as a fully established research journal. With this process in place, a group of mathematicians can now come to the AMR with a proposal, and we can provide a turnkey diamond open access journal, as we have already done for five further journals.

This article explains what we have learned and what we have built. The hope is that others will
join these efforts—but even for those who prefer to build independently, everything here is
documented so that they can do so.

\section{Is There Really a Problem?}

\subsection{Open access versus diamond open access}

``Open access'' means that readers do not pay. The question is: who does?

Under the dominant ``gold'' open access model, costs are shifted to authors through article
processing charges, typically ranging from \$2{,}000 to well over \$10{,}000 per article at major
publishers. These costs are paid from grants, institutional funds, or library ``transformative
agreements''—the same strained budgets that previously struggled under subscription pricing.
Besides this model's potential for incentives other than scientific standards, it
also excludes large segments of the research community:
mathematicians at teaching-focused institutions, independent researchers, early-career faculty
between grants, and many researchers outside North America and Western Europe. A \$3{,}000
article processing charge that is routine at an R1 institution can be prohibitive elsewhere.

Diamond open access takes a different approach: no fees to readers and no fees to authors. When
paired with minimal infrastructure and the service work mathematicians already perform, it is the
only model that is genuinely open to everyone.
During one of my hundred conversations, a colleague asked earnestly: ``But how will you convince libraries to subscribe to your new journals?'' To which I replied: ``Don't you get it? The jig is up!'' The question revealed how deeply the old model has shaped our thinking. Under diamond open access, there are no subscriptions. The libraries are not the customer, because there is no customer. The work is simply \emph{available}.

\subsection{What about the arXiv?}

The arXiv is indispensable for dissemination, and it has dramatically improved access to
mathematical research. But this does not make paywalled journals benign, because the arXiv does
not, on its own, replace the functions that journals serve. Not all papers are posted there, and
even when they are, the versions that mathematicians rely on may differ from the final, refereed
versions. Precise citations to theorem or equation numbers can become unreliable, and
corrections to published work may go unnoticed. Most importantly, the arXiv provides no
organized peer review: journals are not merely dissemination platforms, but mechanisms by
which the community collectively vouches for correctness and quality. The arXiv is a preprint
server, not a validation system.

\section{Highlights of Prior Reform Efforts}

Mathematicians have been discussing the economics of publishing since at least the 1990s. In the
October 2000 \emph{Notices of the AMS}, Allyn Jackson \cite{Jackson2000} listed more than a dozen articles from the
preceding five years addressing pricing, access, and reform. The Cost of Knowledge declaration
in 2012 labeled the situation ``unsustainable and unacceptable.'' 
The diagnosis is not new, and the community has not only talked;
here is but a sampling of concrete efforts that deserve mention.

The Centre Mersenne, hosted by Université Grenoble Alpes and supported by the CNRS, publishes
fourteen diamond open access mathematics journals, including \emph{Algebraic Combinatorics}.
It is a serious and successful operation, supported by substantial government funding.

The Episciences platform, also publicly funded in France, supports overlay journals that conduct
peer review on papers hosted on the arXiv or HAL. It is an elegant model, and one that works well
in practice. It is worth noting that the cost structures of preprint servers and journals are fundamentally different. The arXiv’s operating costs reflect the scale at which it functions, handling an enormous volume of submissions and downloads. A journal publishing even hundreds of papers per year, by contrast, can host all those PDF files at essentially negligible cost.

MathOA, founded in 2016, advocates for and helps fund journal ``flips,'' in which editorial boards
leave commercial publishers to start fair open access alternatives. These efforts have had real
successes, but each flip requires significant organizational and financial effort.

In some subfields, coordinated action has been possible. The Combinatorics Consortium emerged
when several editorial boards resigned from Springer and Elsevier journals to create diamond open
access replacements. This works well in that community, but it relies on particular cultural and
structural features that do not generalize easily.

Finally, publishers such as MSP and the European Mathematical
Society Press have adopted ``Subscribe to Open,'' converting journals to open access once
subscription revenue meets a threshold. This avoids article processing charges and preserves
quality, but it still depends on continued library subscriptions. If those subscriptions falter, access
reverts to being closed.

These are all valuable contributions. But they tend to rely on substantial institutional funding, be
confined to specific subfields, or depend on the same subscription budgets that created the
problem in the first place. The central obstacle is not disagreement about goals, but the absence
of a replicable institutional framework that makes action easy for ordinary editorial boards.

\section{What the AMR Model Looks Like}

\subsection{What is a journal, really?}

At its core, a journal consists of an editorial board of prominent researchers and a commitment to rigorous peer review.
Reputation, trust, and value derive from the judgment of editors and the quality of the papers they
select; everything else—hosting, DOI assignment, copyediting, archiving—is infrastructure. Much
of the traditional publishing apparatus was designed for a print era, and today these functions
can be provided at a fraction of the cost. Indexing and recognition follow from sustained quality
and serious operation over time; they are consequences, not prerequisites.

\subsection{The cost--time tradeoff}

It is possible to run a journal for essentially zero dollars: host it on a university website, manage
submissions by email, and track papers in a spreadsheet. But this does not eliminate costs—it
merely shifts them. In this model, the price is paid in the \emph{time} of the editors, who must perform
large amounts of administrative and clerical work that could otherwise be handled by a proper
editorial management system.

At the other extreme, professionally managed journals with dedicated staff can cost \$100{,}000
per year or more, minimizing the time burden on editors but requiring significant institutional
backing.

The AMR aims for a sweet spot: editors’ time is treated as the most precious resource, and
modest financial investment is used to eliminate administrative overhead, while many traditional
publishing costs are simply cut.

\subsection{The \$3{,}000-per-year stack}

Our current infrastructure costs approximately \$3{,}000 per journal per year. It consists of three
components:
\begin{itemize}
\item EditFlow (from MSP)
 is widely regarded by experienced editors and referees as providing a substantially better editorial workflow than other available systems.
Editors can see the state of the entire journal at a glance,
track refereeing requests and decisions, and manage submissions efficiently without devoting
significant time to administrative bookkeeping;
\item Open Journal Systems (from the Public Knowledge Project) provides a stable,
professional web presence for published papers; and
\item Crossref registers article metadata and assigns persistent DOIs, while ISSNs are obtained through the Library of Congress.
\end{itemize}

Together, these tools replace a large amount of manual coordination and save vast amounts of time, while satisfying the
technical requirements—persistent identifiers, stable hosting, and reliable editorial workflows—of 
a professional journal operation. 

\subsection{Typesetting and format}

At this cost level, journals are fully online; there is no print edition, and professional typesetting
of the kind offered by large publishers is not feasible. Instead, editors perform light copyediting,
with assistance from AI tools for spelling and clarity, after which authors are given direct access
to their \LaTeX\ proofs—typically via collaborative tools such as Overleaf—to finalize the
presentation themselves within the journal’s style guidelines. This is not a degradation of
quality, and in practice it is often an improvement. With the exception of a few publishers—notably
MSP, which employs PhD mathematicians as outstanding typesetters, at commensurate expense—professional production workflows frequently introduce new errors at each stage, requiring
repeated rounds of correction.
(And many of us have had the infuriating experience of working tirelessly to place images \emph{exactly} where we want them to be, getting page breaks in \emph{just} the right place, and so on, only to receive galleys that completely disregard and nullify all our efforts.)
By contrast, our approach allows mathematicians to control exactly
how their work appears, while still ensuring a uniform and professional presentation across the
journal.

\subsection{Editorial and Refereeing Work as Academic Service}

The vast majority of editorial and refereeing work in mathematics is performed by research
mathematicians at universities, who view this as part of their professional service to the
community. Diamond open access does not increase this burden; it simply avoids paying
publishers for the privilege of donating that work.

The only thing worse than not paying for editorial judgment is paying for it.  When revenue depends on acceptances, editorial independence is compromised.
Service-based editorial and refereeing work, for all its limitations, helps keep decision-making
insulated from financial pressure.

\subsection{Editorial Volume and Backlogs}

Traditional journals often impose artificial scarcity: fixed page budgets, strict volume limits, and informal quotas that constrain how many papers can (or must) be published in a given year. These constraints are largely inherited from print-era economics or from prestige models that treat scarcity itself as a signal of quality. The result is familiar to many authors: long backlogs, extended delays, and rejections driven by considerations other than the merit of the work itself.

Under the online-only model, these constraints disappear. There is no need to limit the number of pages published or to manufacture scarcity by arbitrarily rejecting papers that meet the journal’s quality standard. Nor will there ever be any pressure to accept work that does not. Our journals publish on a regular schedule, with each issue containing as many or as few papers as were accepted during the relevant period.

As a consequence, AMR journals have neither backlogs nor quotas. Papers are evaluated individually and accepted or rejected on the basis of quality alone, without regard to how many papers have already been accepted in a given area or time frame. The goal is not rapid publication at the expense of rigor, but the elimination of delays and distortions that arise when editorial decisions are constrained by page budgets or prestige signaling rather than scientific judgment.

\subsection{Continuity and oversight}

A natural concern is sustainability. What happens if an editor-in-chief burns out or an editorial
board disperses?

This is where an overarching organization plays a crucial role. The AMR monitors its journals, maintains
infrastructure, and plans for succession. The probability that both a journal and its parent
organization fail simultaneously is much smaller than either failing alone. Lightweight does not
mean fragile.

\subsection{Funding}

The AMR is not asking authors, readers, or libraries to fund these journals. Because operating costs are so low, the organization has already raised sufficient funds to support its journals for many years, and is in the process of building an endowment to fund them in perpetuity. At \$3{,}000 per year per journal, long-term sustainability is not a serious obstacle.

\section{Why New Institutions?}

The American Mathematical Society (AMS) has served the mathematical community admirably for
over a century and continues to do so. It has also launched a diamond open access journal,
\emph{Communications of the AMS}, which is a welcome and laudable development. However, \emph{CAMS} required
a substantial initial donation and operates at a significantly higher cost structure than the model
described here. Given funding pressures, it seems unlikely that this model can be replicated broadly across the discipline, nor, as discussed above, is it a necessary level of resource expenditure to move diamond open access forward.

More broadly, established organizations accumulate staff, workflows,
and institutional commitments that make certain kinds of change difficult, regardless of
leadership quality. A journal that costs \$3{,}000 per year to run looks very different to an
organization built around legacy production pipelines than to one starting from scratch. This is not a criticism; it is a common feature of paradigm shifts. New models are often easier to
build outside existing structures. The AMR exists not to compete with established societies, but
because certain experiments are easier to conduct independently. If successful, these
experiments may eventually inform larger organizations’ future efforts.

\section{Conclusion: the Experiment is Ongoing}

To date, the AMR has supported diamond open access journals in three distinct ways. First, we
have launched new journals from the ground up, including \emph{JAMR}, the \emph{Journal of Open
Mathematical Problems}, \emph{Applied \& Computational Topology \& Geometry}, and
\emph{Mathematics of Data, Learning, and Intelligence}. Second, we have taken over the
publication of existing journals previously issued by commercial publishers; in particular, the
\emph{Arnold Mathematical Journal}, formerly published by Springer, is now published diamond open access by
the AMR while remaining under the ownership of the Institute for Mathematical Sciences at Stony Brook University. Third, we have facilitated a journal
flip, when the editorial board of \emph{Experimental Mathematics} resigned en masse and founded
the diamond open access \emph{Journal of Experimental Mathematics}
(\url{https://jexpmath.org}).

The AMR itself is an experiment. Whether it succeeds depends on the mathematical community.
Will mathematicians submit strong work to these journals? Will hiring committees, grant panels, and prize committees come to regard publications in these journals as comparably prestigious? These outcomes are not within our
control. We make no claims and offer no prescriptions for how publications should
be weighed. Read the papers. Decide whether the work is good. That is ultimately the only
standard any journal can meet.

Commercial publishers will, of course, continue to exist, and many of the most esteemed
journals in mathematics will remain in their portfolios for the foreseeable future. The AMR’s
response to this is deliberately modest: to do nothing. Rather than
attempting to displace existing journals or to accelerate change through pressure, the aim is
simply to continue building new diamond open access journals and to run them well. If these
journals attract strong editorial boards, publish high-quality work, and are taken seriously by the
community, then over time they can earn the same standing as their subscription-based
counterparts. In this view, time is the primary lever. Institutional norms in mathematics change
slowly, but they do change, and sustained success is ultimately more persuasive than any single
reform effort.

For decades, mathematicians have said that they want their work to be freely available. The
infrastructure now exists to make this possible at minimal cost. The remaining question is not
technical or financial, but collective: whether we are willing to give such institutions the time and social capital
they need to succeed.
\\

\noindent
{\bf Acknowledgements:}
I'm grateful to Joel Hass, Rob Kirby, and Rich Schwartz for many comments and suggestions improving this work.
\\

\noindent
{\bf Conflicts of Interest:}
The author is the Managing Editor of \emph{JAMR}, a founding member of the Board of Directors of the AMR, and its 2026-2027 President.

\bibliographystyle{plain}
\bibliography{AKbibliog}

\end{document}